\documentclass{article}
\usepackage[a4paper, total={5.5in, 8in}]{geometry}
\usepackage{amsmath, amssymb, amsthm}
\usepackage{paralist}
\usepackage[misc]{ifsym}
\usepackage{epsfig} 
\usepackage{epstopdf} 
\usepackage[colorlinks=true]{hyperref}
\hypersetup{urlcolor=blue, citecolor=red}
\allowdisplaybreaks

\usepackage{bbm}
\newtheorem{theorem}{Theorem}[section]

\newtheorem{proposition}[theorem]{Proposition}

\theoremstyle{definition}
\newtheorem{example}[theorem]{Example}
\newtheorem{remark}[theorem]{Remark}

\usepackage{tikz}

\begin{document}
\begin{center}
    {\huge Existence and explicit formula for a~semigroup related to some network problems with unbounded edges}
\end{center}
\medskip
\begin{center}
    Adam B{\l}och\\
    \small{Institute of Mathematics, {\L}{\'o}d{\'z} University of Technology, Poland}\\
    \small{e-mail: adam.bloch@p.lodz.pl}
\end{center}

\medskip

\begin{abstract}
In this paper we consider an initial-boundary value problem related to some network dynamics where the underlying graph has unbounded edges. We show that there exists a $C_0$-semigroup for this problem using a general result from the literature. We also find an explicit formula for this semigroup. This is achieved using the method of characteristics and then showing that the Laplace transform of the solution is equal to the resolvent operator of the generator.
\end{abstract}

\section{Introduction}
In this paper we consider an initial-boundary value problem for a system of linear hyperbolic partial differential equations with some of them posed on bounded spatial intervals and other on unbounded spatial intervals. Let us begin with notation. We introduce three sets of indices: $J_b:=\{1,\ldots,m\}, J_{\infty}^-:=\{1,\ldots, q\}, J_{\infty}^+:=\{1,\ldots r\}$, where $m,q,r\in\mathbb{N}$. Each of them corresponds to a vector-valued function $\boldsymbol{u}(x,t):=\left(u_j(x,t)\right)_{j\in J_b}, \boldsymbol{v}(x,t):=\left(v_j(x,t)\right)_{j\in J_{\infty}^-}, \boldsymbol{w}(x,t):=\left(w_j(x,t)\right)_{j\in J_{\infty}^+}$, where $u_j:[0,1]\to\mathbb{R}$ and $v_j,w_j:[0,\infty)\to\mathbb{R}$. Finally, let $\mathcal{B}$ be a real matrix of dimension $(m+q)\times(m+r)$, which we write as
$$
\mathcal{B}:=\left(\begin{array}{cc}
\mathcal{B}^{11}&\mathcal{B}^{12}\\
\mathcal{B}^{21}&\mathcal{B}^{22}
\end{array}\right),
$$
where $\mathcal{B}^{11}$ is a square block of dimension $m\times m$, $\mathcal{B}^{12}$ is of dimension $m\times r$, $\mathcal{B}^{21}$ of dimension $q\times m$ and $\mathcal{B}^{22}$ of dimension $q\times r$. We consider the following initial-boundary value problem:

\begin{subequations}\label{ibvp}
\begin{align}
    \partial_t\boldsymbol{u}(x,t)&=-\partial_x\boldsymbol{u}(x,t),\quad x\in(0,1), \;t>0,\label{boundedtransport}\\
    \partial_t\boldsymbol{v}(x,t)&=-\partial_x\boldsymbol{v}(x,t),\quad x>0,\; t>0,\label{unboundedtransportminus}\\
    \partial_t\boldsymbol{w}(x,t)&=\partial_x\boldsymbol{w}(x,t),\quad \;x>0, \;t>0,\label{unboundedtransportplus}\\
  \left(\boldsymbol{u}(x,0),\boldsymbol{v}(x,0),\boldsymbol{w}(x,0)\right)^T&=\left(\mathring{\boldsymbol{u}}(x),\mathring{\boldsymbol{v}}(x),\mathring{\boldsymbol{w}}(x)\right)^T,\quad x\in(0,1),\\
    \binom{\boldsymbol{u}(0,t)}{\boldsymbol{v}(0,t)}&=\mathcal{B}\binom{\boldsymbol{u}(1,t)}{\boldsymbol{w}(0,t)},\quad t>0.\label{bc}
\end{align}
\end{subequations}

The paper is organized as follows. Below, we explain how \eqref{ibvp} is related to network problems and comment on the existence of a semigroup solving the abstract Cauchy problem associated with \eqref{ibvp}. In Section \ref{formulasection} we construct the explicit formula for the solution. Section \ref{resolventsection} is devoted to calculating the resolvent operator of the generator of the semigroup. Finally, in Section \ref{theoremsection} we prove that the formula constructed in Section \ref{formulasection} indeed represents the semigroup solution.

\subsection{Relation to network problems} Dynamical systems on networks have been extensively studied since early 1980s. Problem \eqref{ibvp} fits into this framework as well. In the special case $J^+_{\infty}=J^-_{\infty}=\emptyset$ it becomes a transport problem on the interval $[0,1]$, which was extensively studied in \cite{JBPN, JBAFPN, Dornetal, EngelKramar, MKES}, also on infinite graphs in \cite{Dorn}. In this case, the boundary condition \eqref{bc} becomes
\begin{equation}\label{boundedtransportbc}
\boldsymbol{u}(0,t)=\mathcal{B}^{11}\boldsymbol{u}(1,t)
\end{equation}
and it was shown in \cite{JBAF} that $\mathcal{B}^{11}$ has the structure of the adjacency matrix of the line graph of the underlying graph. The equation of the form \eqref{boundedtransport} describes a flow from $0$ to $1$, hence the value of the function at $0$ must be resolved by a boundary condition, while the value at $1$ is completely determined and can be used to describe the values at $0$. This situation is reflected in the boundary condition \eqref{boundedtransportbc}. We also note that the equation \eqref{boundedtransport} is considered on the widely used in the literature interval $[0,1]$ purely for convenience. If needed, each component $u_j$ can be considered on an interval of an arbitrary length and in such a case we can convert the problem to a problem on the interval $[0,1]$ by rescaling.

In \cite{JBABI} the authors consider the so-called telegraph equation on networks. The construction of boundary conditions leads to a special case of a port-Hamiltonian system (see \cite{JacobZwart}) which can be seen as a bidirectional transport problem, i.e., with flows in both directions from $0$ to $1$ and from $1$ to $0$. In the latter case, the value at $1$ must be resolved by a boundary condition while the value at $0$ is determined. As noted in \cite[Section 4]{JBABI}, in an abstract setting a bidirectional transport problem can be converted to unidirectional transport problem by 'inverting' the interval $[0,1]$, which can be seen as a reversing the paramterization of the edges with transport from $1$ to $0$. Precisely speaking, we could define a bounded operator acting as $\upsilon(x) \mapsto \upsilon(1-x)$, which provides similarity relation between corresponding semigroups. We could include transport from $1$ to $0$ in \eqref{ibvp}, however, having in mind the possibility described above, we decided to consider only transport from $0$ to $1$ to keep otherwise cumbersome calculations as simple as possible.

In \cite{EngelKramar} the authors considered graphs contiaining also unbounded edges. In this case, if a transport occurs in the positive direction (towards infinity), then the value at $0$ must be determined, while for the transport in the negative direction (towards $0$) the value at $0$ is determined and provides information at the boundary. This explains the form of the boundary condition \eqref{bc}. We note that contrary to bounded intervals, here we cannot reverse the parametrization to handle only transport in one direction, we need to consider both \eqref{unboundedtransportminus} and \eqref{unboundedtransportplus}. Clearly, in general, the matrix $\mathcal{B}$ need not reflect the structure of a graph but can be an arbitrary matrix of dimension $(m+q)\times (m+r)$, hence network dynamics is only a special class of problems that can be described by \eqref{ibvp}.  

The coefficients appearing next to the spatial derivatives in the equations of the form \eqref{boundedtransport}-\eqref{unboundedtransportplus}, which in our case are all equal to 1, describe the transport velocity. In general, they can be $x$-dependent. This case is considered in \cite{JBABI, EngelKramar}. However, calculating an explicit formula for the semigroup seems to be impossible in this case. In principle, one could use the method of characteristic and iterate the procedure used in the proof of \cite[Theorem 4.2]{JBABI} to construct the solution for any finite time $t$, but the calculations become very involved and still, obtaining a closed-form formula for a solution for any $t\geq 0$ is rather unattainable. On the other hand, it is relatively simple if all of the transport velocities are equal and constant. Then, it is possible to convert the problem to an equivalent problem with unit velocities by simply rescaling the time as $\tau=ct$, where $c$ is the common velocity. There is also an intermediate setup, in which it is assumed that the velocities are linearly dependent over the set $\mathbb{Q}$; we refer to \cite[p. 150]{MKES} and \cite[Eq. (13)]{JBPN} for the constant case and \cite[Eq. (9)]{JBABIII} for the $x$-dependent case. Under this assumption it is again possible to convert the problem to an equivalent problem with unit velocities by introducing artificial edges and vertices. For details we refer to \cite[Section 3.1]{JBABIII}. Since the main result of this paper is the explicit formula, we decided to consider unit velocities, however, we emphasize, that the existence of a semigroup with $x$-dependent velocities can be also handled using the results of \cite{EngelKramar}, as in Section \ref{existencesubsection}.

We also note that the advantage of having an explicit formula for a semigroup is twofold. Firstly, from the engineering point of view it allows to calculate the value of the solution without involving complicated numerical methods. Secondly, from the mathematical point of view it allows to investigate asymptotic behaviour without involving spectral theory, see, for example, \cite{JB, JBABIII, JBPN}.

\subsection{Existence of a semigroup}\label{existencesubsection}
Let us define an operator $(A_p, D(A_p))$ in $\boldsymbol{X}_p:=(L^p(0,1))^{m}\times(L^p(0,\infty))^{q+r}$ by the formula
$$
A_p:=\operatorname{diag}(\underbrace{-\partial_x,\ldots,-\partial_x}_{m+q\;\text{times}},\underbrace{\partial_x,\ldots,\partial_x}_{r\;\text{times}}),
$$
$$
D(A_p):=\left(\left(\begin{array}{c}\boldsymbol{u}\\\boldsymbol{v}\\\boldsymbol{w}\end{array}\right)\in\left(W^{1,p}(0,1)\right)^{m}\times\left(W^{1,p}(0,\infty)\right)^{q+r}:\binom{\boldsymbol{u}(0)}{\boldsymbol{v}(0)}=\mathcal{B}\binom{\boldsymbol{u}(1)}{\boldsymbol{w}(0)}\right)
$$
for any $p\in[1,\infty)$. This operator governs an abstract Cauchy problem related to \eqref{ibvp}. The existence of a semigroup generated by this operator was already proved in \cite{EngelKramar}, however, the authors of this paper consider very general boundary conditions, which make the results on the generation of a semigroup quite involved. Thus, we will show that the existence of a semigroup for \eqref{ibvp} follows from \cite[Corollary 2.18]{EngelKramar}.
\begin{proposition}
    The operator $(A_p,D(A_p))$ generates a strongly continuous semigroup of operators in the space $\boldsymbol{X}_p$ for any $p\in[1,\infty)$.
\end{proposition}
\begin{proof}
We begin with translating the notation of \cite{EngelKramar}. The authors consider $m$ '\emph{external edges}' corresponding to the interval $[0,1]$ and $\ell$ '\emph{internal edges}' corresponding to the interval $[0,\infty)$. Thus, in our case $m$ denotes the same object, while $\ell=q+r$. Next, $f^e, f^i$ are the vector functions defined on the external and internal edges, respectively, thus $f^e=(\boldsymbol{v},\boldsymbol{w}), f^i=\boldsymbol{u}$.

The authors of \cite{EngelKramar} consider an operator of the form
$$
G=c(\cdot)\partial_x,
$$
where $c$ is in general an $x$-dependent, nondiagonal, square matrix of dimension $m+q+r$. In our case this matrix is constant and diagonal, which significantly simplifies the considerations. Precisely speaking,
$$
c(x)\equiv c=\operatorname{diag}(\underbrace{-1,\ldots,-1}_{m+q\;\text{times}},\underbrace{1,\ldots,1}_{r\;\text{times}}).
$$
In particular, the function matrices $q^e, q^i$, which in \cite{EngelKramar} are diagonalizing matrices for $c$, here are the identity matrices of suitable dimensions. Further, $P^i_+, P^i_-$ denote the spectral projections onto the eigenspaces associated with positive and negative eigenvalues, respectively, of the internal part of $c$, that is, the top-left block of $c$ of dimension $m$. Clearly, in our case there are only negative eigenvalues, thus $P^i_+$ vanishes, while $P^i_-$ is again the identity matrix of dimension $m$, as $c$ is already diagonal. Similarly, $P^e_+, P^e_-$ denote the spectral projections onto the eigenspaces associated with positive and negative eigenvalues, respectively, of the external part of $c$, that is, the bottom-right block of $c$ of dimension $q+r$. Thus,
$$
P^e_+=
\left(\begin{array}{cc}
\boldsymbol{0}&\boldsymbol{0}\\
\boldsymbol{0}&\mathcal{I}_r
\end{array}\right),\qquad
P^e_-=
\left(\begin{array}{cc}
\mathcal{I}_q&\boldsymbol{0}\\
\boldsymbol{0}&\boldsymbol{0}
\end{array}\right),
$$
where $\mathcal{I}_q, \mathcal{I}_r$ are the identity matrices of dimension $q$ and $r$, respectively.

The authors of \cite{EngelKramar} consider a 'boundary space' given by $\partial\boldsymbol{X}_p:=\operatorname{rng}(P^e_-)\times\mathbb{R}^m$ and denote
$$
\widetilde{q}:=\operatorname{dim}\partial\boldsymbol{X}_p=\operatorname{rank}(P^e_-)+m
$$
(we added the tilde sign to distinguish it from our $q$). The boundary condition in \cite{EngelKramar} is given in the form
\begin{equation}\label{engkrambc}
V^e_0f^e(0)+V^i_0f^i(0)-V^i_1f^i(1)=Bf,
\end{equation}
where $V^e_0$ is a matrix of dimension $\widetilde{q}+\ell$, $V^i_0$ and $V^i_1$ are matrices of dimension $\widetilde{q}\times m$ and $B$ is a linear and bounded operator from $\boldsymbol{X}_p$ to $\partial\boldsymbol{X}_p$. Clearly, in our case $\operatorname{rank}P^e_i=q$ and hence the boundary space is of dimension $m+q$, which agrees with the number of equations in \eqref{bc}. We need to write the boundary condition in $D(A_p)$ in the form \eqref{engkrambc}. Clearly, $B=0$. Recall that $f^e(0)=(\boldsymbol{v}(0),\boldsymbol{w}(0)), f^i(0)=\boldsymbol{u}(0)$ and $f^i(1)=\boldsymbol{u}(1)$. Hence,
\begin{align*}
\binom{\boldsymbol{u}(0)}{\boldsymbol{v}(0)}-\mathcal{B}\binom{\boldsymbol{u}(1)}{\boldsymbol{w}(0)}&=\binom{\mathcal{I}_{m}}{\boldsymbol{0}}\boldsymbol{u}(0)+\binom{\boldsymbol{0}}{\mathcal{I}_q}\boldsymbol{v}(0)-\binom{\mathcal{B}^{11}}{\mathcal{B}^{21}}\boldsymbol{u}(1)-\binom{\mathcal{B}^{12}}{\mathcal{B}^{22}}\boldsymbol{w}(0)\\
&=\left(\begin{array}{cc}\boldsymbol{0}&-\mathcal{B}^{12}\\
\mathcal{I}_q&-\mathcal{B}^{22}\end{array}\right)\binom{\boldsymbol{v}(0)}{\boldsymbol{w}(0)}+\binom{\mathcal{I}_{m}}{\boldsymbol{0}}\boldsymbol{u}(0)-\binom{\mathcal{B}^{11}}{\mathcal{B}^{21}}\boldsymbol{u}(1),
\end{align*}
and we have
$$
V^e_0=\left(\begin{array}{cc}\boldsymbol{0}&-\mathcal{B}^{12}\\
\mathcal{I}_q&-\mathcal{B}^{22}\end{array}\right),\quad V^i_0=\binom{\mathcal{I}_{m}}{\boldsymbol{0}},\quad V^i_1=\binom{\mathcal{B}^{11}}{\mathcal{B}^{21}}.
$$
We note that the dimension of the first matrix above is $(m+q)\times (q+r)=\widetilde{q}\times\ell$, while of the remaining two $(m+q)\times m=\widetilde{q}\times m$, in agreement with \cite{EngelKramar}. The condition ensuring the existence of a semigroup is that the matrix
$$
\mathcal{R}_0:=\left(V^e_0q^e(0),V^i_1q^i(1)P^i_+-V^i(0)q^i(0)P^i_-\right)
$$
understood as a linear operator in $\partial\boldsymbol{X}_p$ is invertible. We note that this matrix is of dimension $\widetilde{q}\times(\widetilde{q}+\ell)$, thus in this context invertibility must be understood as the existence of a right-inverse, although it is not clearly stated in \cite{EngelKramar}. This, however, agrees with intuition: as a matrix is right-invertible if and only if it has full (row) rank, all of the equations in \eqref{engkrambc} must be linearly independent so that all of the boundary values that must be determined by a boundary condition are uniquely determined. In our case
$$
\mathcal{R}_0=\left(V^e_0,-V^i_0\right)=\left(\begin{array}{ccc}
\boldsymbol{0}&-\mathcal{B}^{12}&-\mathcal{I}_m\\
\mathcal{I}_q&-\mathcal{B}^{22}&\boldsymbol{0}
\end{array}\right)
$$
and clearly the rank of this matrix equals $q+m$ which is equal to its row dimension, thus it has full row rank. We conclude that there exists a $C_0$-semigroup generated by the operator $(A_p, D(A_p))$ for any $p\in[1,\infty)$.
\end{proof}

\begin{figure}[ht]
\begin{center}
\begin{tikzpicture}
    \draw (-2,0) node {$\bullet$};
    \draw (2,0) node {$\bullet$};
    \draw (-2.1,-0.4) node {$\mathbf{v}_1$};
    \draw (2.1,-0.4) node {$\mathbf{v}_2$};
    \draw (0,0.9) node {$u_1$};
    \draw (0,-0.9) node {$u_2$};
    \draw (3.5,-0.3) node {$w_1$};
    \draw (3.5,0.3) node {$v_2$};
    \draw (-3.5,0.2) node {$v_1$};

    \draw [->] (-1.9,0.1) .. controls (0,0.7) .. (1.9,0.1);
    \draw [<-] (-1.9,-0.1) .. controls (0,-0.7) .. (1.9,-0.1);
    \draw [->] (-2.1,0)--(-4.5,0);
    \draw [<-] (2.2,-0.1)--(4.5,-1);
    \draw [->] (2.2,0.1)--(4.5,1);

    \draw (-4.8,0) node {$\cdots$};
    \draw (4.8,-1) node {$\cdots$};
    \draw (4.8,1) node {$\cdots$};
\end{tikzpicture}
\end{center}
\caption{Graph for Example \ref{example}}
\label{graph}
\end{figure}
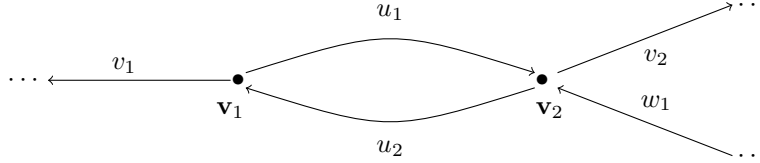

\begin{example}\label{example}
    Let us consider the graph presented in Figure \ref{graph}. The direction of the edges correspond to the direction of the flows. We have $m=q=2$ and $r=1$. At the vertex $\mathbf{v}_1$ we need two boundary conditions to determine the values $u_1(0,t)$ and $v_1(0,t)$ using $u_2(1,t)$. We take them to be
    \begin{align*}
        u_1(0,t)&=\frac{1}{2}u_2(1,t),\\
        v_1(0,t)&=\frac{1}{2}u_2(1,t).
    \end{align*}
    At the vertex $\mathbf{v}_2$ again we need two boundary conditions to determine the values $u_2(0,t)$ and $v_2(0,t)$, but this time we can use two incoming values $u_1(1,t), w_1(0,t)$. Let us take
    \begin{align*}
        u_2(0,t)&=\frac{1}{2}(u_1(1,t)+w_1(0,t)),\\
        v_2(0,t)&=\frac{1}{2}(u_1(1,t)+w_1(0,t)).
    \end{align*}
    Note that if the functions describe density of some physical substances, then these boundary conditions express the Kirchhoff's law and equipartition of the incoming mass between the outgoing edges at each vertex. We have
    \begin{equation}\label{examplebc}
    \left(\begin{array}{c}
    u_1(0,t)\\
    u_2(0,t)\\
    v_1(0,t)\\
    v_2(0,t)
    \end{array}\right)
    =\mathcal{B}\left(\begin{array}{c}
    u_1(1,t)\\
    u_2(1,t)\\
    w_1(0,t)
    \end{array}\right)
    =
    \left(\begin{array}{ccc}
    0&\frac{1}{2}&0\\
    \frac{1}{2}&0&\frac{1}{2}\\
    0&\frac{1}{2}&0\\
    \frac{1}{2}&0&\frac{1}{2}
    \end{array}\right)
    \left(\begin{array}{c}
    u_1(1,t)\\
    u_2(1,t)\\
    w_1(0,t)
    \end{array}\right),
    \end{equation}
    and the blocks of the matrix $\mathcal{B}$ are
    $$
    \mathcal{B}^{11}=\mathcal{B}^{21}=\left(\begin{array}{cc}0&\frac{1}{2}\\\frac{1}{2}&0 \end{array}\right),\quad\mathcal{B}^{12}=\mathcal{B}^{22}=\left(\begin{array}{c}0\\\frac{1}{2} \end{array}\right).
    $$
    Next, we have
    $$
    V^e_0=\left(\begin{array}{ccc}
    0&0&0\\
    0&0&-\frac{1}{2}\\
    1&0&0\\
    0&1&-\frac{1}{2}
    \end{array}\right),\quad V^i_0=\left(\begin{array}{cc}
    1&0\\
    0&1\\
    0&0\\
    0&0
    \end{array}\right),\quad V^i_1=\left(\begin{array}{cc}
    0&\frac{1}{2}\\
    \frac{1}{2}&0\\
    0&\frac{1}{2}\\
    \frac{1}{2}&0
    \end{array}\right),
    $$
    and the boundary condition \eqref{examplebc} written in the form \eqref{engkrambc} is 
    \begin{align*}
    \left(\begin{array}{ccc}
    0&0&0\\
    0&0&-\frac{1}{2}\\
    1&0&0\\
    0&1&-\frac{1}{2}
    \end{array}\right)
    \left(\begin{array}{c}
    v_1(0,t)\\
    v_2(0,t)\\
    w_1(0,t)
    \end{array}\right)
    &+
    \left(\begin{array}{cc}
    1&0\\
    0&1\\
    0&0\\
    0&0
    \end{array}\right)
    \left(\begin{array}{cc}
    u_1(0,t)\\
    u_2(0,t)
    \end{array}\right)\\
    &-
    \left(\begin{array}{cc}
    0&\frac{1}{2}\\
    \frac{1}{2}&0\\
    0&\frac{1}{2}\\
    \frac{1}{2}&0
    \end{array}\right)
    \left(\begin{array}{cc}
    u_1(1,t)\\
    u_2(1,t)
    \end{array}\right)
    =0.
    \end{align*}
    Finally,
    $$
    \mathcal{R}_0=\left(\begin{array}{ccccc}
    0&0&0&-1&0\\
    0&0&-\frac{1}{2}&0&-1\\
    1&0&0&0&0\\
    0&1&-\frac{1}{2}&0&0
    \end{array}\right)
    $$
    and, clearly, $\operatorname{rank}\mathcal{R}_0=4=m+q$.
\end{example}

\section{The explicit formula}\label{formulasection}
To find the solution to the problem \eqref{ibvp}, we use the well-known method of characteristics. Each of the functions $\boldsymbol{u}, \boldsymbol{v}$ and $\boldsymbol{w}$ can be treated separately. We begin with the last function, for which finding the formula is the easiest and which appears in the formulae for the remaining two functions. Each component $w_j$ obeys a simple transport equation with a positive velocity, i.e., in the negative direction (towards 0). Since in this case no boundary condition is required (which is also reflected in \eqref{bc}), the solution is given by the shift of the initial condition:
$$
\boldsymbol{w}(x,t)=\mathring{\boldsymbol{w}}(x+t), \quad x+t>0.
$$

We continue with the formula for $\boldsymbol{u}$. In this case, each component $u_j$ obeys a simple transport equation on the interval $(0,1)$ with negative transport velocity, and thus describes a flow from $0$ to $1$. From the standard theory we know that the solution consists of shifts of the initial condition and the boundary condition depending on the relation between $x$ and $t$. Precisely speaking,
$$
\boldsymbol{u}(x,t)=\mathring{\boldsymbol{u}}(x-t),\quad 0<x-t<1,
$$
and 
$$
\boldsymbol{u}(x,t)=\boldsymbol{\varphi}(t-x), \quad 0<t-x<1,
$$
where $\boldsymbol{\varphi}$ is a function defining the boundary condition $\boldsymbol{u}(0,t)=\boldsymbol{\varphi}(t)$. From \eqref{bc} we have
$$
\boldsymbol{u}(0,t)=\mathcal{B}^{11}\boldsymbol{u}(1,t)+\mathcal{B}^{12}\boldsymbol{w}(0,t).
$$
Note that $\boldsymbol{u}(1,t)=\mathring{\boldsymbol{u}}(1-t)$ and $\boldsymbol{w}(0,t)=\mathring{\boldsymbol{w}}(t)$. Hence, in our case the boundary function $\boldsymbol{\varphi}$ is of the form
$$
\boldsymbol{\varphi}(t)=\mathcal{B}^{11}\mathring{\boldsymbol{u}}(1-t)+\mathcal{B}^{12}\mathring{\boldsymbol{w}}(t)
$$
and thus the solution is given by
$$
\boldsymbol{u}(x,t)=\boldsymbol{\varphi}(t-x)=\mathcal{B}^{11}\mathring{\boldsymbol{u}}(1-t+x)+\mathcal{B}^{12}\mathring{\boldsymbol{w}}(t-x).
$$
This formula is valid as long as the functions $\mathring{\boldsymbol{u}}$ and $\mathring{\boldsymbol{w}}$ are defined, that is, when the conditions $0<1-t+x<1$ and $t-x>0$ are satisfied. However, if the first one is satisfied, then so is the second one, hence, finally, the above formula holds for $0<1-t+x<1$. Using again the boundary condition \eqref{bc} and the formula for $\boldsymbol{u}$ from the first iteration above, we have
\begin{align*}
\boldsymbol{u}(0,t)&=\mathcal{B}^{11}\boldsymbol{u}(1,t)+\mathcal{B}^{12}\boldsymbol{w}(0,t)\\
&=\mathcal{B}^{11}\left(\mathcal{B}^{11}\mathring{\boldsymbol{u}}(2-t)+\mathcal{B}^{12}\mathring{\boldsymbol{w}}(t-1)\right)+\mathcal{B}^{12}\mathring{\boldsymbol{w}}(t)\\
&=\left(\mathcal{B}^{11}\right)^2\mathring{\boldsymbol{u}}(2-t)+\mathcal{B}^{11}\mathcal{B}^{12}\mathring{\boldsymbol{w}}(t-1)+\mathcal{B}^{12}\mathring{\boldsymbol{w}}(t).
\end{align*}
Hence,
using the general solution formula, we have
$$
\boldsymbol{u}(x,t)=\left(\mathcal{B}^{11}\right)^2\mathring{\boldsymbol{u}}(2-t+x)+\mathcal{B}^{11}\mathcal{B}^{12}\mathring{\boldsymbol{w}}(t-x-1)+\mathcal{B}^{12}\mathring{\boldsymbol{w}}(t-x)
$$
and this formula is valid for $0<2-t+x<1$. Iterating the procedure described above, we conclude that the function $u$ is given by the formula
$$
\boldsymbol{u}(x,t)=
\left(\mathcal{B}^{11}\right)^n\mathring{\boldsymbol{u}}(n-t+x)+\sum\limits_{k=0}^{n-1}(\mathcal{B}^{11})^k\mathcal{B}^{12}\mathring{\boldsymbol{w}}(t-x-k),\quad 0<n-t+x<1,
$$
where we understand $\sum_{k=0}^{-1}=0$. Note that the submatrix $\mathcal{B}^{11}$ is square, hence the powers $\left(\mathcal{B}^{11}\right)^k$ are well-defined.

Finally, we consider the formula for $\boldsymbol{v}$. As above, each component again obeys a simple transport equation with a negative velocity with the only difference being the unbounded spatial domain. For small $t<x$ we have
$$
\boldsymbol{v}(x,t)=\mathring{\boldsymbol{v}}(x-t).
$$
For $t>x$, we need to use the boundary condition \eqref{bc}, which gives us
$$
\boldsymbol{v}(0,t)=\mathcal{B}^{21}\boldsymbol{u}(1,t)+\mathcal{B}^{22}\boldsymbol{w}(0,t).
$$
As we can see, $\boldsymbol{v}$ depends only on $\boldsymbol{u}$ and $\boldsymbol{w}$, so contrary to the case of $\boldsymbol{u}$, we do not need to iterate the procedure of obtaining the solution. Instead, we can use the formulae we have already constructed to get
\begin{align*}
\mathcal{B}^{21}\boldsymbol{u}(1,t)+\mathcal{B}^{22}\boldsymbol{w}(0,t)&=\mathcal{B}^{21}\left(\mathcal{B}^{11}\right)^n\mathring{\boldsymbol{u}}(n-t+1)\\
&+\mathcal{B}^{21}\sum\limits_{k=0}^{n-1}\left(\mathcal{B}^{11}\right)^k\mathcal{B}^{12}\mathring{\boldsymbol{w}}(t-1-k)+\mathcal{B}^{22}\mathring{\boldsymbol{w}}(t).
\end{align*}
Hence, from the general formula to the transport equation,
\begin{align*}
\boldsymbol{v}(x,t)&=\mathcal{B}^{21}\left(\mathcal{B}^{11}\right)^n\mathring{\boldsymbol{u}}(n-t+x+1)\\
&+\mathcal{B}^{21}\sum\limits_{k=0}^{n-1}\left(\mathcal{B}^{11}\right)^k\mathcal{B}^{12}\mathring{\boldsymbol{w}}(t-x-k-1)+\mathcal{B}^{22}\mathring{\boldsymbol{w}}(t-x),
\end{align*}
which is valid for $t>x$, where $n$ is chosen so that the condition $0<n-t+x+1<1$ is satisfied. Note that in such a case also $\mathring{\boldsymbol{w}}(t-x-k-1)$ is well defined, i.e., $t-x-k-1>0$, as $k\leq n-1$.

Finally, we conclude that the semigroup $T:=\left(T(t)\right)_{t\geq 0}$ governing the solution to the abstract Cauchy problem associated with \eqref{ibvp} should be given by
\begin{equation}\label{semigroup}
T(t)\left(\begin{array}{c}\mathring{\boldsymbol{u}}\\\mathring{\boldsymbol{v}}\\\mathring{\boldsymbol{w}}\end{array}\right)(x)=\left(\begin{array}{c}\boldsymbol{u}(x,t)\\\boldsymbol{v}(x,t)\\\boldsymbol{w}(x,t)\end{array}\right),
\end{equation}
where particular components are given as follows:
\begin{align*}
    \boldsymbol{u}(x,t)&=\left(\mathcal{B}^{11}\right)^n\mathring{\boldsymbol{u}}(n-t+x)+\sum\limits_{k=0}^{n-1}(\mathcal{B}^{11})^k\mathcal{B}^{12}\mathring{\boldsymbol{w}}(t-x-k),\\
    \boldsymbol{v}(x,t)&=\mathbbm{1}_{\{x>t\}}\mathring{\boldsymbol{v}}(x-t)+\mathbbm{1}_{\{x<t\}}\left(\mathcal{B}^{21}\left(\mathcal{B}^{11}\right)^n\mathring{\boldsymbol{u}}(n-t+x+1)\right.\\
    &\left.+\mathcal{B}^{21}\sum\limits_{k=0}^{n-1}\left(\mathcal{B}^{11}\right)^k\mathcal{B}^{12}\mathring{\boldsymbol{w}}(t-x-k-1)+\mathcal{B}^{22}\mathring{\boldsymbol{w}}(t-x)\right),\\
    \boldsymbol{w}(x,t)&=\mathring{\boldsymbol{w}}(x+t),
\end{align*}
where $0<n-t+x<1$ for $\boldsymbol{u}$ and $0<n-t+x+1<1$ for $\boldsymbol{v}, n\in\mathbb{N}\cup\{0\}$.

\begin{remark}
    We note that the value of $n$ determined by the conditions $0<n-t+x<1$ and $0<n-t+x+1<1$ need not be the same, which follows from the fact that $x\in(0,1)$ for $\boldsymbol{u}$, while $x>0$ for $\boldsymbol{v}$. Its value should be calculated separately depending on the choice of the value of $x$ and $t$. However, since the role played by $n$ in both cases is analogous, we decided not to introduce separate notation.
\end{remark}

\section{The resolvent operator of the generator}\label{resolventsection}
In this section we find the resolvent operator of the generator $(A_p, D(A_p))$ of the semigroup $(T(t))_{t\geq 0}$. In what follows, we use consistent with former considerations boldface notation for vectors, i.e., $\boldsymbol{y}=(y_j)_{j\in J}$ for any set of indices $J$, where the components of the vector and the set of indices are clear from the context.

Let us recall that the resolvent operator, denoted by $R(\lambda, A_p)$, is given by $R(\lambda, A_p):=(\lambda I-A_p)^{-1}$ and is defined for those $\lambda\in\mathbb{C}$, for which the inverse exists and is a bounded operator. Here, $I$ denotes the identity operator in the underlying space. Thus, to calculate $R(\lambda, A_p)$, we need to solve the equation
$$
(\lambda I-A_p)\left(\begin{array}{c}\boldsymbol{u}\\\boldsymbol{v}\\\boldsymbol{w}\end{array}\right)=\left(\begin{array}{c}\boldsymbol{f}\\\boldsymbol{g}\\\boldsymbol{h}\end{array}\right),
$$
where $(\boldsymbol{f},\boldsymbol{g},\boldsymbol{h})^T\in\boldsymbol{X}_p$. Rewriting it in the scalar form, we get the following system of first order linear differential equations:
\begin{align*}
    \lambda u_j(x)+\partial_xu_j(x)&=f_j(x),\quad x\in(0,1), \;j\in J_b,\\
    \lambda v_j(x)+\partial_xv_j(x)&=g_j(x),\quad x\in(0,\infty),\;j\in J_{\infty}^-,\\
    \lambda w_j(x)-\partial_xw_j(x)&=h_j(x),\quad x\in(0,\infty),\;j\in J_{\infty}^+.
\end{align*}
Using the variation of constant formula we have
\begin{align*}
u_j(x)&=e^{-\lambda x}u_0^j+\int\limits_0^xe^{-\lambda(x-s)}f_j(s)\,ds,\\
v_j(x)&=e^{-\lambda x}v_0^j+\int\limits_0^xe^{-\lambda(x-s)}g_j(s)\,ds,
\end{align*}
where $u_0^j, v_0^j$ are certain constants to be determined from the boundary condition appearing in the definition of $D(A_p)$. The equations for $w_j$ are solved by an integrating factor. Multiplying by $e^{-\lambda x}$ and integrating from $x$ to infinity we get
$$
w_j(x)=\int\limits_x^{\infty}e^{\lambda(x-s)}h_j(s)\,ds.
$$
The following considerations are based on the ideas introduced in \cite{JBPN} and extended in \cite{JBAFPN, JBABI}; we refer for details there. To write the solution in an operator form, we denote
$$
\boldsymbol{E}_l(x):=\operatorname{diag}\left(e^x,\ldots,e^x\right),
$$
where $l\in\mathbb{N}$ is the dimension of the (square) matrix. Then,
\begin{align*}
\boldsymbol{u}(x)&=\boldsymbol{E}_m(-\lambda x)\boldsymbol{u}_0+\int\limits_0^x\boldsymbol{E}_m(-\lambda(x-s))\boldsymbol{f}(s)\,ds,\\
\boldsymbol{v}(x)&=\boldsymbol{E}_q(-\lambda x)\boldsymbol{v}_0+\int\limits_0^x\boldsymbol{E}_q(-\lambda(x-s))\boldsymbol{g}(s)\,ds,\\
\boldsymbol{w}(x)&=\int\limits_x^{\infty}\boldsymbol{E}_r(\lambda(x-s))\boldsymbol{h}(s)\,ds.
\end{align*}
Using the boundary condition from the definition of $D(A_p)$, we get
\begin{align*}
\boldsymbol{u}_0&=\mathcal{B}^{11}\boldsymbol{E}_m(-\lambda)\boldsymbol{u}_0+\mathcal{B}^{11}\int\limits_0^1\boldsymbol{E}_m(-\lambda(1-s))\boldsymbol{f}(s)\,ds+\mathcal{B}^{12}\int\limits_0^{\infty}\boldsymbol{E}_r(-\lambda s)\boldsymbol{h}(s)\,ds,\\
\boldsymbol{v}_0&=\mathcal{B}^{21}\boldsymbol{E}_m(-\lambda)\boldsymbol{u}_0+\mathcal{B}^{21}\int\limits_0^1\boldsymbol{E}_m(-\lambda(1-s))\boldsymbol{f}(s)\,ds+\mathcal{B}^{22}\int\limits_0^{\infty}\boldsymbol{E}_r(-\lambda s) \boldsymbol{h}(s)\,ds.
\end{align*}
We rewrite the first equation as
$$
\left(\mathcal{I}_m-\mathcal{B}^{11}\boldsymbol{E}_m(-\lambda)\right)\boldsymbol{u}_0=\mathcal{B}^{11}\int\limits_0^1\boldsymbol{E}_m(-\lambda(1-s))\boldsymbol{f}(s)\,ds+\mathcal{B}^{12}\int\limits_0^{\infty}\boldsymbol{E}_r(-\lambda s)\boldsymbol{h}(s)\,ds,
$$
where $\mathcal{I}_m$ is the identity matrix of size $m$. Taking $\lambda$ with sufficiently large real part, we can decrease the norm of the operator $\mathcal{I}_m-\mathcal{B}^{11}\boldsymbol{E}_m(-\lambda)$ below $1$ to make it invertible with the inverse expressible as the Neumann series. Then,
\begin{align*}
\boldsymbol{u}_0&=\sum\limits_{n=0}^{\infty}\left(\mathcal{B}^{11}\right)^n\boldsymbol{E}_m(-\lambda n)\mathcal{B}^{11}\int\limits_0^1\boldsymbol{E}_m(-\lambda(1-s))\boldsymbol{f}(s)\,ds\\
&+\sum\limits_{n=0}^{\infty}\left(\mathcal{B}^{11}\right)^n\boldsymbol{E}_m(-\lambda n)\mathcal{B}^{12}\int\limits_0^{\infty}\boldsymbol{E}_r(-\lambda s)\boldsymbol{h}(s)\,ds.
\end{align*}
Next, we insert the above formula into the formula for $\boldsymbol{v}_0$ to obtain
\begin{align*}
    \boldsymbol{v}_0&=\mathcal{B}^{21}\boldsymbol{E}_m(-\lambda)\sum\limits_{n=0}^{\infty}\left(\mathcal{B}^{11}\right)^n\boldsymbol{E}_m(-\lambda n)\mathcal{B}^{11}\int\limits_0^1\boldsymbol{E}_m(-\lambda(1-s))\boldsymbol{f}(s)\,ds\\
    &+\mathcal{B}^{21}\boldsymbol{E}_m(-\lambda)\sum\limits_{n=0}^{\infty}\left(\mathcal{B}^{11}\right)^n\boldsymbol{E}_m(-\lambda n)\mathcal{B}^{12}\int\limits_0^{\infty}\boldsymbol{E}_r(-\lambda s)\boldsymbol{h}(s)\,ds\\
    &+\mathcal{B}^{21}\int\limits_0^1\boldsymbol{E}_m(-\lambda(1-s))\boldsymbol{f}(s)\,ds+\mathcal{B}^{22}\int\limits_0^{\infty}\boldsymbol{E}_r(-\lambda s)\boldsymbol{h}(s)\,ds.
\end{align*}
Observe that shifting the index of summation in the first summand from $n+1$ to $n$ and incorporating the third summand into this sum we get
\begin{align*}                      
    \boldsymbol{v}_0&=\mathcal{B}^{21}\sum\limits_{n=0}^{\infty}\left(\mathcal{B}^{11}\right)^n\boldsymbol{E}_m(-\lambda n)\int\limits_0^1\boldsymbol{E}_m(-\lambda(1-s))\boldsymbol{f}(s)\,ds\\
    &+\mathcal{B}^{21}\boldsymbol{E}_m(-\lambda)\sum\limits_{n=0}^{\infty}\left(\mathcal{B}^{11}\right)^n\boldsymbol{E}_m(-\lambda n)\mathcal{B}^{12}\int\limits_0^{\infty}\boldsymbol{E}_r(-\lambda s)\boldsymbol{h}(s)\,ds\\
    &+\mathcal{B}^{22}\int\limits_0^{\infty}\boldsymbol{E}_r(-\lambda s)\boldsymbol{h}(s)\,ds.
\end{align*}
Finally, the resolvent operator $R(\lambda, A_p)$ is given by
\begin{equation}\label{resolvent}
R(\lambda, A_p)\left(\begin{array}{c}\boldsymbol{f}\\\boldsymbol{g}\\\boldsymbol{h}\end{array}\right)(x)=\left(\begin{array}{c}\boldsymbol{u}(x)\\\boldsymbol{v}(x)\\\boldsymbol{w}(x)\end{array}\right),
\end{equation}
where
\begin{align*}
    \boldsymbol{u}(x)&=\boldsymbol{E}_m(-\lambda x)\sum\limits_{n=0}^{\infty}\left(\mathcal{B}^{11}\right)^{n+1}\boldsymbol{E}_m(-\lambda n)\int\limits_0^1\boldsymbol{E}_m(-\lambda(1-s))\boldsymbol{f}(s)\,ds\\
    &+\boldsymbol{E}_m(-\lambda x)\sum\limits_{n=0}^{\infty}\left(\mathcal{B}^{11}\right)^n\boldsymbol{E}_m(-\lambda n)\mathcal{B}^{12}\int\limits_0^{\infty}\boldsymbol{E}_r(-\lambda s)\boldsymbol{h}(s)\,ds\\
    &+\int\limits_0^x\boldsymbol{E}_m(-\lambda(x-s))\boldsymbol{f}(s)\,ds,\\
    \boldsymbol{v}(x)&=\boldsymbol{E}_q(-\lambda x)\mathcal{B}^{21}\sum\limits_{n=0}^{\infty}\left(\mathcal{B}^{11}\right)^n\boldsymbol{E}_m(-\lambda n)\int\limits_0^1\boldsymbol{E}_m(-\lambda(1-s))\boldsymbol{f}(s)\,ds\\
    &+\boldsymbol{E}_q(-\lambda x)\mathcal{B}^{21}\sum\limits_{n=0}^{\infty}\left(\mathcal{B}^{11}\right)^n\boldsymbol{E}_m(-\lambda(n+1))\mathcal{B}^{12}\int\limits_0^{\infty}\boldsymbol{E}_r(-\lambda s)\boldsymbol{h}(s)\,ds\\
    &+\mathcal{B}^{22}\int\limits_0^{\infty}\boldsymbol{E}_r(-\lambda( s+x))\boldsymbol{h}(s)\,ds+\int\limits_0^x\boldsymbol{E}_q(-\lambda(x-s))\boldsymbol{g}(s)\,ds,\\
    \boldsymbol{w}(x)&=\int\limits_x^{\infty}\boldsymbol{E}_r(\lambda(x- s))\boldsymbol{h}(s)\,ds.
\end{align*}

\section{The main theorem}\label{theoremsection}
In this section we prove that the formula obtained in Section \ref{formulasection} indeed represents the semigroup generated by the operator $(A_p,D(A_p))$. Before we start, let us note that since $\boldsymbol{E}(x)$ is a diagonal matrix with equal diagonal entries being exponential functions, it commutes with any other square matrix and $\boldsymbol{E}(x+y)=\boldsymbol{E}(x)\boldsymbol{E}(y)$. We also recall a simple fact from linear algebra that for any matrix $\mathcal{C}$ of dimensions $l_1\times l_2$ there holds
\begin{equation}\label{lemma}
\boldsymbol{E}_{l_1}(x)\mathcal{C}=\mathcal{C}\boldsymbol{E}_{l_2}(x),
\end{equation}
which follows from the fact that $\boldsymbol{E}_l(x)$ is a square matrix with equal diagonal entries.

We are ready to prove the following
\begin{theorem} The Laplace transform $\mathcal{L}\{T\}$ of the family $T=(T(t))_{t\geq 0}$ defined in \eqref{semigroup} satisfies
$$
R(\lambda, A_p)=\mathcal{L}\{T\}(\lambda).
$$
In consequence, $(T(t))_{t\geq 0}$ is the $C_0$-semigroup generated by the operator $(A_p,D(A_p))$.
\end{theorem}

\begin{proof}
Let us recall that, using our notation,
\begin{equation}\label{laplacetransform}
\mathcal{L}\{T\}(\lambda)(x)=\int\limits_0^{\infty}\boldsymbol{E}_{m+q+r}(-\lambda t)T(t)\left(\begin{array}{c}\boldsymbol{f}\\\boldsymbol{g}\\\boldsymbol{h}\end{array}\right)(x)\,dt.
\end{equation}
Similarly as in the former considerations, each component can be treated separately. We denote them by $\left(\mathcal{L}\{T\}\right)_i,\;i=1,2,3$. We begin with with the first one, for which we have
\begin{equation}\label{laptransfcomp1}
\int\limits_0^{\infty}\boldsymbol{E}_m(-\lambda t)\left(\left(\mathcal{B}^{11}\right)^n\boldsymbol{f}(n-t+x)+\sum\limits_{k=0}^{n-1}(\mathcal{B}^{11})^k\mathcal{B}^{12}\boldsymbol{h}(t-x-k)\right)\,dt.
\end{equation}
The first summand is treated analogously as in the proof of \cite[Proposition 3.3]{Dorn}: 
\begin{align*}
&\int\limits_0^{\infty}\boldsymbol{E}_m(-\lambda t)\left(\mathcal{B}^{11}\right)^n\boldsymbol{f}(n-t+x)\,dt\\
=&\int\limits_0^x\boldsymbol{E}_m(-\lambda t)\boldsymbol{f}(x-t)\,dt
+\sum\limits_{n=1}^{\infty}\;\int\limits_{n+x-1}^{n+x}\boldsymbol{E}_m(-\lambda t)\left(\mathcal{B}^{11}\right)^n\boldsymbol{f}(n-t+x)\,dt\\
=&\int_0^x\boldsymbol{E}_m(-\lambda(x-s))\boldsymbol{f}(s)\,ds\\
+&\boldsymbol{E}_m(-\lambda x)\sum\limits_{n=0}^{\infty}\left(\mathcal{B}^{11}\right)^{n+1}\boldsymbol{E}_m(-\lambda n)\int\limits_0^1\boldsymbol{E}_m(-\lambda(1-s))\boldsymbol{f}(s)\,ds.
\end{align*}
We only note that the limits of integration in the sum of integrals come from solving the condition $0<n-t+x<1$ with respect to $t$. In the second summand in \eqref{laptransfcomp1} we split the integral into the sum of integrals as above, then change the order of summation and substitute $s=t-x-k$:
\begin{align*}
&\int\limits_0^{\infty}\boldsymbol{E}_m(-\lambda t)\sum\limits_{k=0}^{n-1}(\mathcal{B}^{11})^k\mathcal{B}^{12}\boldsymbol{h}(t-x-k)\,dt\\
=&\sum\limits_{n=0}^{\infty}\sum\limits_{k=0}^{n-1}\;\int\limits_{n+x-1}^{n+x}\boldsymbol{E}_m(-\lambda t)\left(\mathcal{B}^{11}\right)^k\mathcal{B}^{12}\boldsymbol{h}(t-x-k)\,dt\\
=&\sum\limits_{k=0}^{\infty}\sum\limits_{n=k+1}^{\infty}\;\int\limits_{n+x-1}^{n+x}\boldsymbol{E}_m(-\lambda t)\left(\mathcal{B}^{11}\right)^k\mathcal{B}^{12}\boldsymbol{h}(t-x-k)\,dt\\
=&\sum\limits_{k=0}^{\infty}\sum\limits_{n=k+1}^{\infty}\;\int\limits_{n-k-1}^{n-k}\boldsymbol{E}_m(-\lambda(s+x+k))\left(\mathcal{B}^{11}\right)^k\mathcal{B}^{12}\boldsymbol{h}(s
)\,dt\\
=&\boldsymbol{E}_m(-\lambda x)\sum\limits_{k=0}^{\infty}\boldsymbol{E}_m(-\lambda k)\left(\mathcal{B}^{11}\right)^k\mathcal{B}^{12}\sum\limits_{n=k+1}^{\infty}\;\int\limits_{n-k-1}^{n-k}\boldsymbol{E}_r(-\lambda s)\boldsymbol{h}(s)\,ds\\
=&\boldsymbol{E}_m(-\lambda x)\sum\limits_{k=0}^{\infty}\boldsymbol{E}_m(-\lambda k)\left(\mathcal{B}^{11}\right)^k\mathcal{B}^{12}\int\limits_{0}^{\infty}\boldsymbol{E}_r(-\lambda s)\boldsymbol{h}(s)\,ds,
\end{align*}
where we used \eqref{lemma} and the obvious fact that for any $k\in\mathbb{N}\cup\{0\}$,
\begin{equation}\label{union}
\bigcup\limits_{n=k+1}^{\infty}[n-k-1,n-k]=[0,\infty).
\end{equation}
Hence, the first component of the Laplace integral \eqref{laplacetransform} is given by
\begin{align*}
    \left(\mathcal{L}\{T\}\right)_1(\lambda)(x)&=\boldsymbol{E}_m(-\lambda x)\sum\limits_{n=0}^{\infty}\left(\mathcal{B}^{11}\right)^{n+1}\boldsymbol{E}_m(-\lambda n)\int\limits_0^1\boldsymbol{E}_m(-\lambda(1-s))\boldsymbol{f}(s)\,ds\\
    &+\boldsymbol{E}_m(-\lambda x)\sum\limits_{n=0}^{\infty}\boldsymbol{E}_m(-\lambda n)\left(\mathcal{B}^{11}\right)^n\mathcal{B}^{12}\int\limits_{0}^{\infty}\boldsymbol{E}_r(-\lambda s)\boldsymbol{h}(s)\,ds\\
    &+\int_0^x\boldsymbol{E}_m(-\lambda(x-s))\boldsymbol{f}(s)\,ds.
\end{align*}
We proceed with the second component, i.e.,
\begin{align*}
\int\limits_0^{\infty}\boldsymbol{E}_q(-\lambda t)&\left(\mathbbm{1}_{\{x>t\}}\boldsymbol{g}(x-t)+\mathbbm{1}_{\{x<t\}}\left(\mathcal{B}^{21}\left(\mathcal{B}^{11}\right)^n\boldsymbol{f}(n-t+x+1)\right.\right)\\
&\left.+\mathcal{B}^{21}\sum\limits_{k=0}^{n-1}\left(\mathcal{B}^{11}\right)^k\mathcal{B}^{12}\boldsymbol{h}(t-x-k-1)+\mathcal{B}^{22}\boldsymbol{h}(t-x)\right)\,dt.
\end{align*}
For convenience, we calculate each of the four summands separately. Substituting $s=x-t$ in the first one and $s=t-x$ in the last one and using \eqref{lemma} we get
$$
\int\limits_0^{\infty}\boldsymbol{E}_q(-\lambda t)\mathbbm{1}_{\{x>t\}}\boldsymbol{g}(x-t)\,dt=\int\limits_0^x\boldsymbol{E}_q(-\lambda(x-s))g(s)\,ds
$$
and
$$
\int\limits_0^{\infty}\boldsymbol{E}_q(-\lambda t)\mathbbm{1}_{\{x<t\}}\mathcal{B}^{22}\boldsymbol{h}(t-x)\,dt=\mathcal{B}^{22}\int\limits_0^{\infty}\boldsymbol{E}_r(-\lambda(s+x))\boldsymbol{h}(s)\,ds.
$$
The second summand is treated analogously as the first summand of \eqref{laptransfcomp1}, however, since the expression in the argument of the function $\boldsymbol{f}$ is slightly different, we need to adjust the limits of integration in the sum. Again using \eqref{lemma}, we have
\begin{align*}
    &\int\limits_0^{\infty}\boldsymbol{E}_q(-\lambda t)\mathbbm{1}_{\{x<t\}}\mathcal{B}^{21}\left(\mathcal{B}^{11}\right)^n\boldsymbol{f}(n-t+x+1)\,dt\\
    =&\sum\limits_{n=0}^{\infty}\;\int\limits_{n+x}^{n+x+1}\boldsymbol{E}_q(-\lambda t)\mathcal{B}^{21}\left(\mathcal{B}^{11}\right)^n\boldsymbol{f}(n-t+x+1)\,dt\\
    =&\boldsymbol{E}_q(-\lambda x)\mathcal{B}^{21}\sum\limits_{n=0}^{\infty}\left(\mathcal{B}^{11}\right)^n\boldsymbol{E}_m(-\lambda n)\int\limits_0^1\boldsymbol{E}_m(-\lambda(1-s))\boldsymbol{f}(s)\,ds.
\end{align*}
The third summand is handled as follows:
\begin{align*}
&\int\limits_0^{\infty}\boldsymbol{E}_q(-\lambda t)\mathbbm{1}_{\{x<t\}}\mathcal{B}^{21}\sum\limits_{k=0}^{n-1}\left(\mathcal{B}^{11}\right)^k\mathcal{B}^{12}\boldsymbol{h}(t-x-k-1)\,dt\\
=&\sum\limits_{n=0}^{\infty}\sum\limits_{k=0}^{n-1}\int\limits_{n+x}^{n+x+1}\boldsymbol{E}_q(-\lambda t)\mathcal{B}^{21}\left(\mathcal{B}^{11}\right)^k\mathcal{B}^{12}\boldsymbol{h}(t-x-k-1)\,dt\\
&=\boldsymbol{E}_q(-\lambda x)\sum\limits_{k=0}^{\infty}\boldsymbol{E}_q(-\lambda(k+1))\sum\limits_{n=k+1}^{\infty}\;\int\limits_{n-k-1}^{n-k}\boldsymbol{E}_q(-\lambda s)\mathcal{B}^{21}\left(\mathcal{B}^{11}\right)^k\mathcal{B}^{12}\boldsymbol{h}(s)\,ds\\
=&\boldsymbol{E}_q(-\lambda x)\mathcal{B}^{21}\sum\limits_{k=0}^{\infty}\left(\mathcal{B}^{11}\right)^k\boldsymbol{E}_m(-\lambda(k+1))\mathcal{B}^{12}\int\limits_{0}^{\infty}\boldsymbol{E}_r(-\lambda s)\boldsymbol{h}(s)\,ds,
\end{align*}
where once again we used \eqref{lemma} and \eqref{union}. Combining the above formulae, the second component of the Laplace integral \eqref{laplacetransform} is given by
\begin{align*}
\left(\mathcal{L}\{T\}\right)_2(\lambda)(x)&=\boldsymbol{E}_q(-\lambda x)\mathcal{B}^{21}\sum\limits_{n=0}^{\infty}\left(\mathcal{B}^{11}\right)^n\boldsymbol{E}_m(-\lambda n)\int\limits_0^1\boldsymbol{E}_m(-\lambda(1-s))\boldsymbol{f}(s)\,ds\\
&+\boldsymbol{E}_q(-\lambda x)\mathcal{B}^{21}\sum\limits_{k=0}^{\infty}\left(\mathcal{B}^{11}\right)^k\boldsymbol{E}_m(-\lambda(k+1))\mathcal{B}^{12}\int\limits_{0}^{\infty}\boldsymbol{E}_r(-\lambda s)\boldsymbol{h}(s)\,ds\\
&+\mathcal{B}^{22}\int\limits_0^{\infty}\boldsymbol{E}_r(-\lambda(s+x))\boldsymbol{h}(s)\,ds+\int\limits_0^x\boldsymbol{E}_q(-\lambda(x-s))\boldsymbol{g}(s)\,ds.
\end{align*}
Finally, calculating the last component is straightforward:
\begin{align*}
\left(\mathcal{L}\{T\}\right)_3(\lambda)(x)&=\int\limits_0^{\infty}\boldsymbol{E}_r(-\lambda t)\boldsymbol{h}(x+t)\,dt\\
&=\int\limits_x^{\infty}\boldsymbol{E}_r(-\lambda(x-s))\boldsymbol{h}(s)\,ds.
\end{align*}
Comparing $\left(\mathcal{L}\{T\}\right)_i, i=1,2,3$, with the formulae for the components $\boldsymbol{u}, \boldsymbol{v}, \boldsymbol{w}$ in \eqref{resolvent} and using a standard result in semigroup theory, \cite[Theorem 1.10]{EngelNagel}, we obtain the assertion.
\end{proof}

\end{document}